\documentclass[a4paper,12pt]{article}

\usepackage{amsthm,amsmath,amscd,mathtools}
\usepackage[vmargin=2cm,hmargin=2cm,headheight=14.5pt,top=2cm,headsep=.5cm]{geometry}
\usepackage{hyperref}
\usepackage[utf8]{inputenc}
\usepackage[autobold]{mathfixs}
\usepackage{stackrel}
\usepackage{cases}
\usepackage{circuitikz}
\usepackage{tikz}

\usepackage[charter]{mathdesign}
\usepackage{bm}
\usepackage{caption}
\usepackage{subcaption}
\usepackage{xcolor}
\usepackage{graphicx}
\DeclareMathSizes{12}{12}{8}{6}

\usepackage{cite}
\usepackage{enumitem}
\setlist[itemize,2]{label=$\centerdot$}
\setlist[itemize,3]{label=$\triangle$}

\newtheoremstyle{ptheorem}{1em}{0em}{\itshape}{}{\bfseries}{.}{.5em}{\thmname{#1}\thmnumber{
 #2}\thmnote{ (\hspace{-.01pt}{#3})}}

\theoremstyle{ptheorem}

\newtheorem{thm}{Theorem}[section]
\newtheorem{pro}[thm]{Proposition}
\newtheorem{lem}[thm]{Lemma}
\newtheorem{cor}[thm]{Corollary}

\newtheoremstyle{hdef}{1em}{0em}{}{}{\bfseries}{.}{.5em}{\thmname{#1}\thmnumber{
 #2}\thmnote{ (\hspace{-.01pt}{#3})}}
\theoremstyle{hdef}

\newtheorem{dfn}[thm]{Definition}
\newtheorem{rem}[thm]{Remark}

\newtheorem{exa}[thm]{Example}

\numberwithin{equation}{section}
\numberwithin{figure}{section}

\allowdisplaybreaks

\usepackage{hyperref}
\begin{document}
 \title{Stieltjes differential systems with non monotonic derivators}

\author{
Marl\`ene Frigon \footnote{This work was partially supported by NSERC Canada.}\\
\normalsize e-mail: frigon@dms.umontreal.ca\\
\normalsize \emph{D\'epartement de math\'ematiques et de statistique,} \\ \normalsize\emph{Universit\'e de Montr\'eal, Canada.}\\
F. Adri\'an F. Tojo\footnote{The author was partially supported by Ministerio de Econom\'ia y Competitividad, Spain, and FEDER, project MTM2013-43014-P, and by the Agencia Estatal de Investigaci\'on (AEI) of Spain under grant MTM2016-75140-P, co-financed by the European Community fund FEDER.} \\
\normalsize e-mail: fernandoadrian.fernandez@usc.es\\
\normalsize \emph{Instituto de Ma\-te\-m\'a\-ti\-cas, Facultade de Matem\'aticas,} \\ \normalsize\emph{Universidade de Santiago de Com\-pos\-te\-la, Spain.}\\
}
\date{}
 \maketitle

\begin{abstract}
In this work we study Stieltjes differential systems of which the derivators are allowed to change sign. This leads to the definition of the notion of \emph{function of controlled variation}, a characterization of precompact sets of $g$-continuous functions, and an explicit expression of $g$-exponential maps. Finally, we prove a Peano-type existence result and apply it to a model of fluid stratification on buoyant miscible jets and plumes.
\end{abstract}


 \section{Introduction}

In the last years, there has been quite an interest in the theory regarding problems with Stieltjes derivatives. This kind of problems can be traced back to the works of Kurzweil~\cite{Kurzweil1957} and later of Schwabik and Tvrd\'y~\cite{Sch1,Sch2,Tvrdy} but, whereas these works deal with Stieltjes measures directly through the associated integral problems, some others, such as~\cite{FP2016,PoMa,PoRo}, center their attention in the definition and the meaning of the Stieltjes derivatives.

In a recent article of {López Pouso and M\'arquez Alb\'es \cite{PoMa2}}, the authors study the case of systems with Stieltjes derivatives with several different derivators. All of them are assumed to be left continuous and nondecreasing. They establish the existence of a solution of the following system of differential equations with respect to $h$:
\begin{equation}\label{hdifeq-1}
x'_h(t)=f(t,x(t))\ \ \mu_h\text{-a.e.}\ t\in[t_0,t_0+\tau],\quad x(0)=x_0,
\end{equation}
where $h=(h_1,\dots,h_n) : [t_0,t_0+T] \to {\mathbb R}^n$ is such that $h_i$ is left continuous and nondecreasing for every $i=1,\dots,n$, (see~\cite[Theorem 4.5]{PoMa2}).

In this manuscript, we consider more general maps $h$. More precisely, we want to see if the nondecreasing assumption on all $h_i$ can be dropped.
In particular, for $g:[t_0,t_0+\tau]\to{\mathbb R}$ a left continuous function of bounded variation, we study the problem
\begin{equation}\label{gdifeq-1}
x'_g(t)=f(t,x(t))\ \ \mu_g\text{-a.e.}\ t\in[t_0,t_0+\tau],\quad x(0)=x_0.
\end{equation}

 First, let us observe that it is much easier to establish the existence of a solution of the weaker integral equation:
 \begin{equation}\label{ginteq}
 x(t)= x_0 + \int_a^tf(s,x(s))\operatorname{d}\mu_g(s),
\end{equation}
where $\mu_g$ is the Lebesgue-Stieltjes measure with respect to $g$.
Indeed, it is well known that there exist two left continuous nondecreasing functions $g_+, g_-$ such that $g=g_+-g_-$. This remark leads us to consider the following system:
\begin{equation}\label{gdifsys}
\begin{aligned}
y'_{g_+}(t) &= f(t,y(t)+z(t)),\\
z'_{g_-}(t) &= -f(t,y(t)+z(t)),\\
(y(t_0),z(t_0)) &= (x_0,0).
\end{aligned}
\end{equation}
For $h=(g_+,g_-)$, {\cite[Theorem~4.5]{PoMa2}, ensures} the existence of a solution $(y,z)$ of~\eqref{gdifsys}. It is easy to see that $(y,z)$ is also a solution of the integral system:
\begin{equation}\label{gintsys}
\begin{aligned}
y(t) &= x_0 +\int_a^tf(t,y(t)+z(t))\operatorname{d}\mu_{g_+},\\
z(t) &= -\int_a^tf(t,y(t)+z(t))\operatorname{d}\mu_{g_-}.
\end{aligned}
\end{equation}
So, $y+z$ satisfies
\[
(y+z)(t)= x_0+ \int_a^tf(t,y(t)+z(t))\operatorname{d}(\mu_{g_+}-\mu_{g_-}),
\]
that is, $x=y+z$ is a solution of~\eqref{ginteq}.

It is important to realize that, although the previous procedure is simple, it is lacking in that it does not provide any information on the form or behavior of the $g$-derivative, and that is because we only obtained a solution of the integral equation~\eqref{ginteq}. What is more, in order to study the problem~\eqref{gdifeq-1} we would first have to define the $g$-derivative, the signed measure $\mu_g$ and also the positive measure $|\mu_g|$ in this context. Of course, the Radon-Nikodym derivative will be properly defined as it follows from the integral equation, but it lacks an interpretation in terms of $g$. In this sense, there have been previous treatments of derivatives with respect to functions of bounded variation from the measure perspective (see~\cite{Daniell,Garg}) or related to the Kurzweil-Stieltjes integral~\cite{Giselle}. Measures have the great advantage of allowing for a general treatment without worrying about the behavior at specific points, but this same asset turns into a disadvantage when it is precisely that trait that we want to control in practical applications.

It is for the reasons mentioned above that we will study the case of systems with $g$-derivators that can change sign, focusing on the obtaining of solutions of a system of $g$-differential equations. In order to achieve this goal, we will have to rework the theory regarding the Fundamental Theorem of Calculus for the Lebesgue-Stieltjes integral as presented in~\cite{PoRo}, where it was done only for measures associated to nondecreasing functions. Allowing our derivators to be non monotonic will force us to introduce new concepts, such as that of function of controlled variation, in order to correctly define a new notion of absolute continuity which is not based on related measures (see~\cite{Garg}) but relies directly on the function.

This development of the theory will allow us to obtain an explicit formula for the $g$-exponential map which will be the solution of a linear $g$-differential equation. Then, we will present an analog of the Peano theorem in this context. Finally, in the last section, we show how the theory of $g$-differential equations with non monotone $g$ can be applied to a model of fluid stratification on buoyant miscible jets and plumes presented in~\cite{Camassa}.

\section{Preliminaries}

We recall some facts about bounded variation functions and their relation with signed measures.

\begin{dfn} \label{df:var-g}
Let $I\subset{\mathbb R}$ be an interval and $g:I\to{\mathbb R}$. Given an interval $J\subset I$, we denote
\begin{equation*}
{\mathcal P}(J): = \{P=(x_1,\dots, x_{k_P})\ :\ k_P\geqslant 2;\ x_j\in J,\ j=1,\dots,k_P;\ x_j\leqslant x_{j+1},\ j=1,\dots,k_P-1\}.
 \end{equation*}
The \emph{total variation of $g$ in $J$} is defined by
 \[
 \operatorname{var}_gJ:=\sup_{P\in\mathcal{P}(J)} \sum_{j=1}^{k_P-1} | g(x_{j+1})-g(x_j) |\in[0,\infty].
 \]
If $\operatorname{var}_g(I)<\infty$, we say that $g$ is a function of \emph{bounded variation} and we denote
\begin{align*}
\operatorname{BV}(I,{\mathbb R}) &= \{g : I \to {\mathbb R} : g \text{ is a function of bounded variation}\},\\
\operatorname{BV}_-(I,{\mathbb R}) &= \{g \in \operatorname{BV}(I,{\mathbb R}) : g \text{ is left continuous}\}.
 \end{align*}
 \end{dfn}

{From now on we write $g(x^+):=\lim_{y\to x^+}g(y)$ and $g(x^-):=\lim_{y\to x^-}g(y)$. In order to avoid having to separate cases when $x$ is an extremity of $I$, for instance, in the case $I=[a,b]$, we define $g(a^-):=g(a)$ and $g(b^+):=g(b)$.}
 It is well known that functions of bounded variation are \emph{regulated}, so lateral limits always exist. For $g \in \operatorname{BV}(I,{\mathbb R})$, let us define
\begin{align*}
D_g^+ &:= \{x\in {I}: g(x^+)-g(x^-)>0\},\\
D_g^- &:= \{x\in {I}: g(x^+)-g(x^-)<0\},\\
D_g &:= D_g^+\cup D_g^-,
\end{align*}
and
\begin{equation}\label{cdelta}
C_g:= \{x\in {I}:\ g(x)=g(y), y\in(x-\varepsilon,x+\varepsilon)\mbox{ for some } \varepsilon\in{\mathbb R}^+\}.
\end{equation}
It is also well known that $D_g$ is countable~\cite{PoRo}. Given an interval $J\subset I$, we consider
 \begin{align*}\operatorname{var}_g^+J=&\sup_{P\in\mathcal{P}(J)} \sum_{j=1}^{k_P-1} [ g(x_{j+1})-g(x_j) ]^+,\\
 \operatorname{var}_g^-J=&\sup_{P\in\mathcal{P}(J)} \sum_{j=1}^{k_P-1} [ g(x_{j+1})-g(x_j) ]^-,
 \end{align*}
where $c^+=\max\{c,0\}$ and $c^-=-\min\{c,0\}$.
 We call $\operatorname{var}_g^+J$ and $\operatorname{var}_g^-J$ the \emph{positive variation} and the \emph{negative variation} of $g$ on $J$ respectively.

 For $g \in BV_-(I,{\mathbb R})$, $\operatorname{var}_g^+$, $\operatorname{var}_g^-$ and $\operatorname{var}_g$ can serve to construct measures $\mu_g^+$, $\mu_g^-$ and $|\mu_g|$ respectively~\cite[Theorem 4.1.9]{Ben} in the following way. Define, for instance, $v_g(x):=\operatorname{var}_g[a,x)$, $x\in(a,b]$. Thus defined, $v_g$ is left continuous and there is a unique measure, $|\mu_g|$ such that
 \[
 |\mu_g|([x,y))=v_g(y)-v_g(x)=\operatorname{var}_g[a,y)-\operatorname{var}_g[a,x)=\operatorname{var}_g[x,y).
 \]
 Observe that
 \[
 |\mu_g|(\{x\})=v_g(x^+)-v_g(x)=:\operatorname{var}_g[x,x^+]= g(x^+)-g(x).
 \]
 So, $|\mu_g|(\{x\}) \ne 0$ if and only if $x \in D_g$.
 Note that $C_g$ is, by definition, an open set in the usual topology of $I$. Therefore, it can be rewritten uniquely as the disjoint countable union of open intervals, say $C_g=\bigcup_{n\in{\mathbb N}} (a_n,b_n)$. Thus, it is easy to see that $|\mu_{ g}|(C_g) = 0$. The measures $\mu_g^+$ and $\mu_g^-$ are defined similarly from $\operatorname{var}_g^-$ and $\operatorname{var}_g^+$ respectively.

 \begin{lem} The measure $|\mu|_g$ satisfies $|\mu_g|=\mu_g^++\mu_g^-$.
\end{lem}

\begin{dfn}[Jordan Decomposition] We define the signed measure $\mu_g:=\mu_g^+-\mu_g^-$.
 \end{dfn}

The definitions presented here for $\mu_g^+$, $\mu_g^-$ and $|\mu_g|$ coincide with the ones stated in order to obtain the Jordan decomposition~\cite[Definition 5.1.2]{Ben}, that is, for any $\mu_g$-measurable set $A$, we have that
\begin{align*}
\mu_g^+(A)=&\sup\{\mu_g(B)\ :\ B\subset A,\ B\ \mu_g\text{-measurable}\},\\
\mu_g^-(A)=& \sup\{-\mu_g(B)\ :\ B\subset A,\ B\ \mu_g\text{-measurable}\}, \\
|\mu_g|(A)=&\ \mu_g^+(A) + \mu_g^-(A).
\end{align*}

We denote by $\operatorname{L}^1_g(I,{\mathbb R})$ the space of $|\mu_g|$-integrable functions on the interval $I$, that is, the set of \hbox{$|\mu_g|$-measurable} maps $h : I \to {\mathbb R}$ such that $\int_I |h|\operatorname{d}|\mu_g|<\infty$. Notice that
\[
\int_I h\operatorname{d}|\mu_g|=\int_{I\backslash D_g} h\operatorname{d}|\mu_g|+\sum_{t\in D_g}h(t)|g(t^+)-g(t)|.
\]

\begin{rem} Observe that if a function is $|\mu_g|$-integrable, it is also $\mu_g$-integrable.
\end{rem}

\section{Functions of controlled variation}

In this section, we introduce the notion of functions of controlled variation.

\begin{dfn}\label{df:CV}
 We say that $g\in\operatorname{BV}_-(I,{\mathbb R})$ is a \emph{function of controlled variation} if there exists a closed set $F_g\subset I$ such that $g$ is monotonic on each connected component of $I\backslash F_g$ and {$|\mu_{ g}|(F_g\backslash D_g)=0$}. We denote
 \[
 \operatorname{CV}_-(I,{\mathbb R}) = \{g \in \operatorname{BV}_-(I,{\mathbb R}) : g \text{ is a function of controlled variation}\}.
\]
\end{dfn}

\begin{pro}
The set	$\operatorname{CV}_-(I,{\mathbb R})$ is a dense subset of $\operatorname{BV}_-(I,{\mathbb R})$.
\end{pro}

\begin{proof}
Rational polynomials are contained in $\operatorname{CV}_-(I,{\mathbb R})$ and are dense in $\operatorname{L^1}(I,{\mathbb R})$. There is a continuous inclusion of $\operatorname{BV}(I,{\mathbb R})$ into $\operatorname{L^1}(I,{\mathbb R})$~\cite{Kreuzer2013}, so $\operatorname{CV}_-(I,{\mathbb R})$ is dense in $\operatorname{BV}_-(I,{\mathbb R})$.
\end{proof}

\begin{rem}
Not all functions in $\operatorname{BV}_-(I,{\mathbb R})$ are of controlled variation. In~\cite[Section 3]{Amo}, they build a continuous function $g$ of bounded variation that is of monotonic type on no interval. That implies that $F_g$ would have to be taken equal to $I$, but in that case $|\mu_g|(F_g\backslash D_g)=|\mu_g|(I)\ne 0$, so $g\not\in\operatorname{CV}_-(I,{\mathbb R})$. Given
$g \in \operatorname{BV}_-(I,{\mathbb R})$, we have that $g=g_1-g_2$, where $g_1$ and $g_2$ are continuous nondecreasing functions, and therefore $g_1, g_2\in\operatorname{CV}_-(I,{\mathbb R})$. Hence, we deduce that $\operatorname{CV}_-(I,{\mathbb R})$ is not a vector space.
Such functions in $\operatorname{BV}_-(I,{\mathbb R})\backslash\operatorname{CV}_-(I,{\mathbb R})$ do not occur in practical applications and are only of interest as counterexamples.
\end{rem}

In Remark~\ref{remwcv}, we will stress why it is important for $g$ to be of controlled variation.

\begin{rem}\label{remhd}
Observe that, for $g \in \operatorname{CV}_-(I,{\mathbb R})$, the role of $F_g$ in the previous definition could be played by $F_g\backslash C_g$, which is compact and such that $(F_g\backslash C_g)\backslash D_g$ is hereditarily disconnected (that is, all of its connected components consist of just one point~\cite[p. 223]{hart}). To see this, take a connected component $E$ of $F_g\backslash D_g$ with more than a point. The connected components in the usual topology of $I$ are intervals and, since $|\mu_g|(E)=0$, $g$ has to be constant on $E$. Therefore, $\mathring E$ can be subtracted from $F_g$ because $g$ is monotonic on $E$ and $F_g\backslash\mathring E$ is closed. The set $E\backslash\mathring E$ has at most two points, each one is a connected component of $(F_g\backslash D_g)\backslash\mathring E$. We can do this process with all connected components of $F_g\backslash D_g$ with more than one point. We arrive to a new $F_g$, which is precisely $F_g\backslash C_g$ and $(F_g\backslash C_g)\backslash D_g$ is hereditarily disconnected. Moreover, $|\mu_g|((F_g\backslash C_g)\backslash D_g)=|\mu_g|(F_g\backslash( D_g\cup C_g))= 0$.
\end{rem}

\begin{rem} Hereditarily disconnected sets, such as the Cantor set, can be uncountable.
 \end{rem}

In what follows, for $g \in \operatorname{CV}_-(I,{\mathbb R})$, we will assume that its associated closed set $F_g$ is a subset of $I\backslash C_g$.

\begin{lem}\label{lemst}
If $g\in \operatorname{CV}_-(I,{\mathbb R})$, then there exists a countable family of pairwise disjoint subintervals of $I$, $\{I_\lambda\}_{\lambda\in\Lambda}$, such that $F_g= I\backslash\bigcup_{\lambda\in\Lambda} I_\lambda$ and, for each $\lambda\in\Lambda$, $g$ is monotonic on $I_\lambda = (a_\lambda,b_\lambda)$.
\end{lem}

\begin{proof}
We assume, using Remark~\ref{remhd}, that $F_g\backslash D_g$ is hereditarily
disconnected. Since $F_g$ is closed, the open set ${\mathbb R}\backslash F_g$, is a (unique) disjoint union of a countable (because the usual topology of ${\mathbb R}$ is second countable) family of pairwise disjoint open intervals, $\bigcup_{\lambda\in\Lambda}I_\lambda$ where $I_\lambda=(a_\lambda,b_\lambda)$ for $\lambda\in\Lambda$. Since each $(a_\lambda,b_\lambda)$ is connected and contained in $I\backslash F_g$, $g$ is monotonic on $(a_\lambda,b_\lambda)$.
 \end{proof}

\begin{dfn}\label{df:Idecomp}
Let $g \in \operatorname{CV}_-(I,{\mathbb R})$ and its associated closed set $F_g= I\backslash\bigcup_{\lambda\in\Lambda}I_\lambda\subset I\backslash C_g$ with $\{I_\lambda\}_{\lambda\in\Lambda}$ a countable family of pairwise disjoint open intervals such that, for each $\lambda\in\Lambda$, $g$ is monotonic on $I_\lambda$. We define
\begin{align*}
\Lambda^+ &= \{\lambda \in \Lambda : g \text{ is nondecreasing on } I_\lambda = (a_\lambda,b_\lambda)\}, \quad& \Lambda^- &= \Lambda\backslash\Lambda^+,\\
A^+ &= \bigcup_{\lambda^+\in\Lambda^+}I_{\lambda^+}, &
A^- &= \bigcup_{\lambda^-\in\Lambda^-}I_{\lambda^-}.
\end{align*}
\end{dfn}

Observe that $
\mu_g^+(A^-) = \mu_g^+(D_g^-) = 0$ and $\mu_g^-(A^+) = \mu_g^-(D_g^+) = 0
$.

Let us recall the existence of a Hahn decomposition associated to $g$.

\begin{thm}[{Hahn Decomposition, \cite[Theorem 5.1.6]{Ben}}] There exists a $\mu_g$-measurable set, $P \subset I$, such that $\mu_g^+(E)=\mu_g(E\cap P)$ and $\mu_g^-(E)=-\mu_g(E\backslash P)$ for every $E\subset I$ $\mu_g$-measurable.
\end{thm}

The sets given in Definition~\ref{df:Idecomp} are related to the Hahn decomposition associated to $\mu_g$.

\begin{pro} Let $g \in \operatorname{CV}_-(I,{\mathbb R})$ and $P \subset I$ the $\mu_g$-measurable set given by the Hahn decomposition of $\mu_g$. We have, for every $E\subset I$ $\mu_g$-measurable,
\[
\mu_g^+(E)=\mu_g(E\cap P) = \mu_g(E\cap(A_g^+\cup D_g^+)) \quad\text{and}\quad \mu_g^-(E)=-\mu_g(E\backslash P)= -\mu_g((E\cap(A_g^-\cup D_g^-)).
\]
\end{pro}

\begin{proof} Let $E\subset I$ $\mu_g$-measurable. Since $|\mu_g|(F_g\backslash D_g)=0$,
\begin{align*}
\mu_g(E\cap P) &= \mu_g^+(E)\\
 &= \mu_g^+(E\cap (F_g\backslash D_g)) + \mu_g^+(E\cap (A_g^+\cup D_g^+)) + \mu_g^+(E\cap (A_g^-\cup D_g^-))\\
 &= \mu_g^+(E\cap (A_g^+\cup D_g^+)).
\end{align*}
Similarly, we deduce that $\mu_g^-(E)=-\mu_g(E\backslash P)= -\mu_g((E\cap(A_g^-\cup D_g^-))$.
\end{proof}

\section{The $g$-continuity and the relative compactness in $BC_g(I,\mathbb R)$}

From now on and unless stated otherwise, we will consider $g\in\operatorname{CV}_-(I,{\mathbb R})$ and $F_g \subset I\backslash C_g$ its associated closed set.

As notation, consider $[x,y]:=[y,x]$ if $y\leqslant x$. If we define $\Delta_g(x,y)=\operatorname{var}_g[x,y]$ for $x,y\in I$, we find that $\Delta$ is a displacement.

\begin{dfn}[\cite{MaTo}]\label{deltadef} Let $X\ne\emptyset$ be a set. A \emph{displacement} is a function $\Delta:X^2\to{\mathbb R}$ such that the following properties hold:
	\begin{itemize}
		\item[(H1)] $\Delta(x,x)=0,\ x\in X$.
		\item[(H2)] For all $x,y\in X$,
		\[
 \stackrel[{z\rightharpoonup y}]{}{\underline\lim}|\Delta(x,z)|=|\Delta(x,y)|,
 \]
		where
		\[ \stackrel[{z\rightharpoonup y}]{}{\underline\lim}|\Delta(x,z)|:=\sup \left\{\liminf_{n\to\infty}|\Delta(x,z_n)|: (z_n)_{n\in{\mathbb N}}\subset X,\ |\Delta(y,z_n)|\xrightarrow{n\to\infty}0\right\},
\]
		and the limit is considered in the usual topology of ${\mathbb R}$.
	\end{itemize}
\end{dfn}

The fact that $\Delta_g$ is a displacement implies that we already know a lot regarding the topology generated by $g$.

\begin{dfn}
	Given $t\in I$ and $r\in{\mathbb R}^+$, we define the \textit{$g$-open ball of center $t$ and radius $r$}, as \hbox{$B_g(t,r):=\{s\in I\ :\ \operatorname{var}_g{[t,s]}<r\}$}.	Also, we define the \emph{$g$-topology} in the following way:
	\[
\tau_g:=\left\{U\subset I\ :\ \forall t\in U\ \exists r\in{\mathbb R}^+,\ B_g(t,r)\subset U\right\}.
\]
\end{dfn}

\begin{rem}This definition coincides with~\cite[Definition 3.1]{FP2016} when $g$ is nondecreasing.
\end{rem}

When necessary, we will denote by $\tau_u$ the usual topology on ${\mathbb R}^n$.

\begin{dfn}\label{cont}
A map $h:I\to {\mathbb R}$ is said to be \emph{$g$-continuous} if $h:(I,\tau_{g})\to({\mathbb R},\tau_u)$ is continuous. We denote by ${\mathcal C}_g(I,{\mathbb R})$ the set of \hbox{$g$-continuous} functions and by $BC_g(I,{\mathbb R})$ the set of bounded $g$-continuous functions.
\end{dfn}

The proof of the following result is analogous to~\cite[Theorem 3.4]{FP2016}.

\begin{pro}
The space $BC_g(I,{\mathbb R})$ endowed with the supremum norm is a Banach space.
\end{pro}

As usual, $g$-continuity can be characterized using $g$-open balls.

\begin{lem}[\cite{MaTo}]\label{lcont}
A map $h:I\to {\mathbb R}$ is $g$-continuous if and only if, for every $t\in I$ and $\varepsilon\in{\mathbb R}^+$, there exists $\delta\in{\mathbb R}^+$ satisfying $|h(s)-h(t)|<\varepsilon$ for every $s\in B_g(t,\delta)$.
\end{lem}

{From now on $I$ will be a compact interval $[a,b]$ in $\mathbb R$.}

	The next proposition can be proven as in~\cite[Proposition 3.2]{FP2016}.

\begin{pro} \label{lc}
If $g \in \operatorname{CV}_-({I},{\mathbb R})$ and $h:I\to{\mathbb R}$ is $g$-continuous, then
	\begin{enumerate}
		 \item[(1)] $h$ is continuous from the left at every $t\in(a,b]$;
		 \item[(2)] if $g$ is continuous at $t\in[a, b)$, then so is $h$;
		 \item[(3)] if $g$ is constant on some $[\alpha,\beta]\subset I$, then so is $h$.
\end{enumerate}
\end{pro}

In~\cite{FP2016}, the authors gave sufficient conditions on certain sets of regulated functions in $BC_g(I,{\mathbb R})$ for the case where $g$ is nondecreasing to ensure that they are relatively compact. Here, we provide a characterization for $BC_g(I,{\mathbb R})$ in our more general setting using Vala's Theorem~\cite{Vala}.

\begin{dfn} Let $X$ be a set and $(Y,d)$ be a metric space. We say that a map $f:X\to Y$ is \emph{precompact} if $f(X)$ is precompact. Let $K(X,Y)$ be the family of precompact maps from $X$ to $Y$ endowed with the usual metric induced by $d$.
\end{dfn}

\begin{thm}[{Ascoli-Arzel\`a-Vala \cite[Theorem 1]{Vala}}]\label{thmvala}
 Let $K(X,Y)$ be endowed with the topology induced by the metric and let $S\subset K(X,Y)$. Then $S$ is precompact if and only if
 \begin{enumerate}
 \item[(i)] $S(x):=\{f(x)\in Y\ :\ f\in S\}$ is precompact for all $x\in X$;
 \item[(ii)] $S$ is of \emph{equal variation}, i.e., for every $\varepsilon\in{\mathbb R}^+$, there is a finite covering $\{X_k\}_{k=1}^m$ of $X$ such that $d(f(x),f(y))<\varepsilon$ for every $f\in S$, every $x,y\in X_k$ and every $k=1,\dots,m$.
 \end{enumerate}
\end{thm}

Now, all we have to do is to translate this characterization to the space $BC_g(I,{\mathbb R})$. To this end, first observe that, since ${\mathbb R}$ is a complete space, a subset $A\subset{\mathbb R}$ is precompact if and only if it is relatively compact~\cite[(3.17.5)]{Die}. Furthermore, ${\mathbb R}$ satisfies Heine-Borel's property, so $A$ is precompact if and only if it is bounded. This implies that
\[
BC_g(I,{\mathbb R})=BC((I,\tau_g),({\mathbb R},\tau_u)){\subset} K((I,\tau_g),({\mathbb R},|\cdot|)),
\]
and, therefore, Theorem~\ref{thmvala} is applicable to $BC_g(I,{\mathbb R})$.

\begin{dfn} A set $S\subset BC_g(I,{\mathbb R})$ is said to be \emph{$g$-equicontinuous} if, for every $\varepsilon\in{\mathbb R}^+$ and $t\in  {I}$, there exists $\delta\in\mathbb{R}$ such that $|f(t)-f(s)|<\varepsilon$ for every $f\in S$ and every $s\in I$ such that $\operatorname{var}_g[t,s]<\delta$. We say that $S$ is \emph{uniformly $g$-equicontinuous} if, for every $\varepsilon\in{\mathbb R}^+$, there exists $\delta\in{\mathbb R}^+$ such that $|f(t)-f(s)|<\varepsilon$ for every $f\in S$ and every $t,s\in I$ such that $\operatorname{var}_g[t,s]<\delta$.
\end{dfn}

\begin{dfn}\label{dgs}
A set $S\subset BC_g(I,{\mathbb R})$ is said to be \emph{$g$-stable} if, for every $t\in[a,b)\cap D_g$ and every $\varepsilon\in{\mathbb R}^+$, there exist $\delta\in{\mathbb R}^+$ and a finite covering $\{X_k\}_{k=1}^m$ of $[t,t+\delta)\cap I$ such that $|f(x)-f(y)|<\varepsilon$ for every $f\in S$, every $x,y\in X_k$ and every $k=1,\dots,m$.
\end{dfn}

\begin{rem}\label{remu} Observe that in Definition~\ref{dgs} we can change $[t,t+\delta)\cap I$ by $(t,t+\delta)\cap I$. Indeed, if $\{X_k\}_{k=1}^m$ is a suitable finite covering of $(t,t+\delta)\cap I$, then $\{t\}\cup \{X_k\}_{k=1}^m$ is an appropriate finite covering of $[t,t+\delta)\cap I$.
\end{rem}

Stability is a property that serves to compensate the lack of compactness of $(I,\tau_g)$. The space $(I,\tau_g)$ is not compact in general because, as it was pointed out in~\cite[Example 3.3]{FP2016}, there are $g$-continuous functions defined on $I$ which are not bounded, something that would always be the case if $(I,\tau_g)$ was compact. Stability has been present in the literature in several ways. In~\cite[Theorem~4.2(iii)]{FP2016}, for regulated functions, they ask $S$ to have uniform right-hand side limits in $t\in[a,b)\cap D_g$, something which coincides with our definition in that case. It also appears in~\cite[Theorem 1]{prz} (named \emph{stability} as well) or~\cite[Section 2.12, p. 62]{Corduneanu} (named \emph{equiconvergence}) with a different formulation consistent with the topology of ${\mathbb R}$.

%
%

\begin{lem} \label{unif-eq=stab}
Let $S\subset BC_g(I,{\mathbb R})$. If $S$ is uniformly $g$-equicontinuous,
 then $S$ is $g$-equicontinuous and $g$-stable.
\end{lem}

\begin{proof} Clearly, $S$ is $g$-equicontinuous. Take $t\in[a,b)\cap D_g$ and fix $\varepsilon\in{\mathbb R}^+$. Let $\rho\in{\mathbb R}^+$ be such that $|f(x)-f(y)|<\varepsilon$ for every $f\in S$ and every $x, y\in I$ such that $\operatorname{var}_g[x,y]<\rho$. Now, there exists $\delta\in{\mathbb R}^+$ such that $\operatorname{var}_g(t,t+\delta)<\rho$ because $g$ is regulated. Thus, by Remark~\ref{remu}, $S$ is $g$-stable.
\end{proof}

\begin{thm}\label{thmbc}
 Let $S\subset BC_g(I,{\mathbb R})$, then $S$ is precompact if and only if
 \begin{enumerate}
 \item[(i)] $S(t)$ is bounded for all $t\in I$;
 \item[(ii)] $S$ is $g$-equicontinuous;
 \item[(iii)] $S$ is $g$-stable.
 \end{enumerate}
\end{thm}

\begin{proof}
 As we have argued before, $S(t)$ is precompact if and only if $S(t)$ is bounded, so Conditions~(i) in both Theorem~\ref{thmvala} and Theorem~\ref{thmbc} coincide.

 Now, let us see that, if (i)--(iii) hold, then $S$ satisfies Condition~(ii) of Theorem~\ref{thmvala}.
 Fix $\varepsilon\in{\mathbb R}^+$. For every $t\in I$, there exists $\delta_t\in{\mathbb R}^+$ such that $|f(t)-f(s)|<\frac{\varepsilon}{2}$ for every $f\in S$ and every $s\in I$ such that $\operatorname{var}_g[t,s]<\delta_t$. Consider $U_t:= B_g(t,\delta_t)$ for every $t\in I$. If $t\in I\backslash D_g$, then, by the continuity of $g$ at $t$, there exists $\rho_t\in{\mathbb R}^+$ such that $V_t=(t-\rho_t,t+\rho_t)\cap I\subset U_t$. On the other hand, if $t\in D_g$, since $g$ is left continuous, there exists $\sigma_t\in{\mathbb R}^+$ such that $(t-\sigma_t,t]\cap I\subset U_t\cap[a,t]$. Finally, since $S$ is $g$-stable, there exists $\eta_t\in{\mathbb R}^+$ and a finite covering $\{X_{t,k}\}_{k=1}^{m_t}$ of $[t,t+\eta_t)\cap I$ such that $|f(x)-f(y)|<\varepsilon$ for every $f\in S$, every $x,y\in X_{t,k}$ and every $k=1,\dots,m_t$. For $t \in D_g$, we set $\rho_t:=\min\{\sigma_t,\eta_t\}$. We denote
 $V_t:=(t-\rho_t,t+\rho_t)\cap I$.
 The set $\{V_t\}_{t\in I}$ is an open cover of $(I,\tau_u)$, so there exists a finite subcover $\{V_{t_i}\}_{i=1}^j$.

 We consider the finite covering of $I$,
 \begin{align*}
 \mathcal{Y} &:= \big\{V_{t_i} : i \in \{1,\dots,j\} \text{ and } t_i \not\in D_g\big\} \cup \big\{V_{t_i}\backslash [t_i,t_i+\rho_{t_i}) : i \in \{1,\dots,j\} \text{ and } t_i \in D_g\big\} \\
 &\quad\cup \big\{[t_i,t_i+\rho_{t_i}) \cap X_{t_i,k_i} : i \in \{1,\dots,j\}, t_i \in D_g \text{ and } k_i \in \{i,\dots,m_{t_i}\}\big\}.
 \end{align*}
 If $x, y \in V_{t_i}$ for some $t_i\not\in D_g$ or if $x, y \in V_{t_i}\backslash [t_i,t_i+\rho_{t_i})$ for some $t_i \in D_g$, then $x, y \in U_{t_i}$ and
 $$
 |f(x)-f(y)|\leqslant|f(x)-f(t_i)|+|f(t_i)-f(y)|<\varepsilon \quad \text{for every $f \in S$}.
 $$
 Also, if $x, y \in [t_i,t_i+\rho_{t_i})\cap X_{t_i,k_i}$ for some $t_i\in D_g$ and some $k_i \in \{1,\dots,m_{t_i}\}$, then $|f(x)-f(y)| < \varepsilon$ for every $f \in S$.
 Therefore, $S$ satisfies~(ii) of Theorem~\ref{thmvala} and hence, $S$ is precompact.

 Now, let us assume that $S$ is precompact. It was noticed before that $S$ satisfies~(i). Let $\varepsilon > 0$. For $t \not\in D_g$ and for $\rho > 0$, we consider the set
 $$
 S_{\varepsilon,t}(\rho):=\{f\in S\ :\ |f(t)-f(s)|<\varepsilon \ \text{for all $s \in I$ such that $\operatorname{var}_g[t,s] < \rho$}\}.
 $$
 It is easy to see that $\{S_{\varepsilon,t}(\rho)\}_{\rho\in{\mathbb R}^+}$ is an open cover of $\overline S$. Since $S$ is precompact, and $BC_g(I,{\mathbb R})$ is complete, $\overline S$ is compact, so there is a finite subcover $\{S_{\varepsilon,t}(\rho_k)\}_{k=1,\dots,m}$. Take $\delta:=\min\{\rho_k : k=1,\dots,m\}$. Then, if $s\in I$ is such that $\operatorname{var}_g[t,s] <\delta$, one has $|f(t)-f(s)|<\varepsilon$ for every $f \in S$. Similarly, for $t \in D_g$, using the left continuity of $f \in S$ at $t$, one can see that $\{S_{\varepsilon,t}(\rho)\}_{\rho\in (0,g(t^+)-g(t))}$ is an open cover of $\overline S$. Again, the precompactness of $S$ {ensures} the existence of a finite subcover $\{S_{\varepsilon,t}(\rho_k)\}_{k=1,\dots,m}$. Thus, if $s\in I$ is such that $\operatorname{var}_g[t,s] < \min\{\rho_k : k=1,\dots,m\}$, one has $|f(t)-f(s)|<\varepsilon$ for every $f \in S$. So, $S$ is $g$-equicontinuous. Finally, Theorem~\ref{thmvala} implies that there exists a finite covering $\{X_k\}_{k=1}^m$ of $I$ such that $|f(x)-f(y)| < \varepsilon$ for every $f \in S$, every $x, y \in X_k$ and every $k \in \{1,\dots,m\}$. Therefore $S$ is $g$-stable.
\end{proof}

Combining Lemma~\ref{unif-eq=stab} and Theorem~\ref{thmbc}, we deduce the following result.

\begin{cor}\label{cor:precompact}
 Let $S\subset BC_g(I,{\mathbb R})$. If
 \begin{enumerate}
 \item[(i)] $S(t)$ is bounded for all $t\in I$ and
 \item[(ii)] $S$ is uniformly $g$-equicontinuous,
 \end{enumerate}
 then $S$ is precompact.
\end{cor}

\section{The $g$-derivative}

We now generalize to non monotonic derivator $g$ the notion of $g$-derivative introduced in~\cite{PoRo}.

\begin{dfn}\label{deriv}
 The \emph{derivative with respect to $g$} (or \emph{$g$-derivative}) of a function $f: I\to {\mathbb R}$ at a point $x\in D_g \cup (I\backslash(F_g\cup\overline{C_g}))$ is defined as follows, provided that the corresponding limits exist:
 \begin{align*}
 f'_g(x) =\begin{dcases}
 \lim_{y \to x}\frac{f(y)-f(x)}{g(y)-g(x)} &\text{if $x\not\in D_{g}$,}\\ \lim_{y \to x^+}\frac{f(y)-f(x)}{g(y)-g(x)} &\text{if $x\in D_{g}$.}\end{dcases}
 \end{align*}
\end{dfn}


\begin{pro}\label{procarder}
Let $f : I \to {\mathbb R}$ and $x_0\in D_g \cup (I\backslash(F_g\cup\overline{C_g}))$. Then the following statements are equivalent:
\begin{enumerate}
\item[(1)] The map $f$ is $g$-differentiable at $x_0$.
\item[(2)] There exist $\alpha \in {\mathbb R}$, a set $N\subset I$ containing $x_0$ and a function $\varphi:N\to{\mathbb R}$ such that $\varphi(x_0)=0$, $g(x)\ne g(x_0)$ for every $x \in N\backslash\{x_0\}$,
 $$f(x)=f(x_0)+[\varphi(x)+\alpha][g(x)-g(x_0)] \quad \text{for every $x \in N$,}
 $$
 and one of the following holds:
 \begin{enumerate}
 \item[(i)] $x_0 \not\in D_g$, $N$ is open in $I$ and $\varphi$ is continuous at $x_0$;
 \item[(ii)] $x_0\in D_g$, $N$ is open in $[x_0,\infty)\cap I$ and $\varphi$ is right continuous at $x_0$.
 \end{enumerate}
\end{enumerate}
In the case~(2), we have that $\alpha=f'_g(x_0)$.
\end{pro}

\begin{proof} $(1)\Rightarrow (2)$: If $f$ is $g$-differentiable at $x_0 \not\in D_g$ (resp. $x_0 \in D_g$), then $g(x)\ne g(x_0)$ in a neighborhood $N$ of $x_0$ in $I$ (resp. in $I\cap [x_0,\infty)$) (otherwise we could not take the limit that defines $f'_g(x_0)$). Hence, it is enough to define
\[
\varphi(x):=\frac{f(x)-f(x_0)}{g(x)-g(x_0)}-f'_g(x_0)
 \]
 for $x\in N\backslash\{x_0\}$ and $\varphi(x_0)=0$. By the definition of the derivative, $\varphi$ is continuous at $x_0$ (resp. right continuous at $x_0$).

$(2)\Rightarrow (1)$: Since $\varphi$ is right continuous at $x_0$,
\[
\lim_{x\to x_0^+}\frac{f(x)-f(x_0)}{g(x)-g(x_0)}=\lim_{x\to x_0^+}(\varphi(x)+\alpha)=\alpha.
\]
In the case where $x_0\not\in D_g$, we have the analogous result with the left limit. So, in either case, there exists $f'_g(x_0)=\alpha$.
\end{proof}

\begin{rem}\label{remwcv} Now, we see why it is important that $g$ is of controlled variation. Assume, for instance, that we have a point $t_0\in (a,b)\backslash D_g$ such that there are two sequences $\{t_n\}, \{s_n\}$ in $(a,t_0)$ converging to $t_0$ and satisfying that $g(t_n)<g(t_0) < g(s_n)$; that is, that $t_0$ does not belong to any open interval where $g$ is monotonic. Furthermore, assume that $g$ is differentiable at $t_0$. Then, if $f$ is $g$-differentiable at $t_0$, we have that
\[
f'(t_0)=\lim_{t\to t_0}\frac{f(t)-f(t_0)}{t-t_0}=\lim_{t\to t_0}\frac{f(t)-f(t_0)}{g(t)-g(t_0)}\frac{g(t)-g(t_0)}{t-t_0}=f'_g(t_0)g'(t_0).
\]
But also, observe that, for every $n\in{\mathbb N}$,
\[
\frac{g(t_n)-g(t_0)}{t_n-t_0}>0 > \frac{g(s_n)-g(t_0)}{s_n-t_0},
\]
which implies that $g'(t_0)=0$. Hence, $f'(t_0)=0$.

In conclusion, the fact that $t_0$ does not belong to any open interval where $g$ is monotonic forces very restrictive conditions on the behavior of $f$.
\end{rem}

The next result establishes a more general chain rule than the one in~\cite{PoRo}.

\begin{pro}[Chain rule]\label{cr} Let $g_1 \in \operatorname{CV}_-\big((a,b),{\mathbb R}\big)$, $g_2 \in \operatorname{CV}_-\big((c,d),{\mathbb R}\big)$, $f: (a,b) \to (c,d)$ and $h : (c,d) \to {\mathbb R}$. Assume that $g_2\circ f$ is $g_1$-differentiable at some $x_0 \in (a,b)$ with $f(x_0) \not\in D_{g_2}$, and $h$ is $g_2$-differentiable at $f(x_0)$. Then, $h\circ f$ is $g_1$-differentiable at $x_0$ and
\[ (h\circ f)_{g_1}'(x_0)=h'_{g_2}(f(x_0))(g_2\circ f)'_{g_1}(x_0).\]
\end{pro}

\begin{proof} We write the proof for $x_0\not\in D_{g_1}$. It is analogous for $x_0\in D_{g_1}$. Observe that, since $(g_2\circ f)_{g_1}'$ exists, $g_1(x)\ne g_1(x_0)$ in a neighborhood of $x_0$.
Since $h$ is $g_2$-differentiable at $f(x_0)\not\in D_{g_2}$, by Proposition~\ref{procarder}, there exist $N$ a neighborhood of $f(x_0)$ and a function $\varphi: N\to{\mathbb R}$ continuous at $f(x_0)$ and such that $\varphi(f(x_0))=0$ and
\[
h(y)=h(f(x_0))+[\varphi(y)+h_{g_2}'(f(x_0))][g_2(y)-g_2(f(x_0))].
\]
Therefore,
\[
\frac{h(f(x))-h(f(x_0))}{g_1(x)-g_1(x_0)}= (\varphi(f(x))+h_{g_2}'(f(x_0)))\frac{g_2(f(x))-g_2(f(x_0))}{g_1(x)-g_1(x_0)}.
\]
Taking the limit, we have that
\[
(h\circ f)_{g_1}'(x_0)=h_{g_2}'(f(x_0))(g_2\circ f)_{g_1}'(x_0).
\]
\end{proof}


The following corollary is obtained by taking $g_2$ as the identity.

\begin{cor}\label{ccr} Let $g \in CV_-(I,{\mathbb R})$, $f:I \to (c,d)$ a function $g$-differentiable at $x_0$ and $h : (c,d) \to {\mathbb R}$ a function differentiable at $f(x_0)$. Then, $(h\circ f)$ is $g$-differentiable at $x_0$ and
$(h\circ f)_g'(x_0)=h'(f(x_0))f'_g(x_0)$.
\end{cor}


\section{The Fundamental Theorem of Calculus}

%

From now on, we consider $g\in\operatorname{CV}_-(I,{\mathbb R})$ and the sets associated to $g$ provided by Definition~\ref{df:Idecomp},
\[
A_g^+ = \bigcup_{\lambda^+\in\Lambda^+}I_{\lambda^+} \quad\text{and}\quad A_g^- = \bigcup_{\lambda^-\in\Lambda^-}I_{\lambda^-}.
\]

\begin{dfn}\label{ac2}
A function $h:I\to{\mathbb R}$ is said to be \emph{$g$-absolutely continuous} if it is $g$-continuous and if, for every $\varepsilon\in{\mathbb R}^+$, there exists $\delta\in{\mathbb R}^+$ such that, for any family $\{(\alpha_\xi,\beta_\xi)\}_{\xi\in\Xi^+}$ of pairwise disjoint subintervals of $A_g^+$, any family $\{(\alpha_\xi,\beta_\xi)\}_{\xi\in\Xi^-}$ of pairwise disjoint subintervals of $A_g^-$, and any sets $E^+\subset D_g^+$, $E^-\subset D_g^-$, we have that, if
\[
\sum_{\xi\in\Xi^+}[g(\beta_\xi)-g(\alpha_\xi^+)] + \sum_{\xi\in\Xi^-}[g(\alpha_\xi^+)-g(\beta_\xi)]
+ \sum_{t\in E^+}[g(t^+)-g(t)] + \sum_{t\in E^-}[g(t)-g(t^+)] <\delta,
\]
 then
 \[
 \sum_{\xi\in\Xi^+\cup \Xi^-}|h(\beta_\xi)-h(\alpha_\xi^+)| + \sum_{t\in E^+\cup E^-}|h(t^+)-h(t)| <\varepsilon.
 \]
 We denote by $\operatorname{AC}_g(I,{\mathbb R})$ the set of $g$-absolutely continuous functions.
\end{dfn}

In the next two theorems, we show that, from $\mu_g$-integrable maps, we can defined $g$-absolutely continuous functions which are $g$-differentiable $|\mu_g|$-almost every where.

\begin{thm}\label{thmint}
Let $v\in\operatorname{L}^1_{g}([a,b))$ and $h(x):=\int_{[a,x)}v(s)\operatorname{d}\mu_g$ \ for $x\in I$. Then $h\in\operatorname{AC}_g(I,{\mathbb R})$.
\end{thm}

\begin{proof} We consider the case $v\geqslant 0$ $|\mu_g|$-a.e. on $I$ since $v$ has to be the
difference of two such functions $|\mu_g|$-a.e. Let $\varepsilon\in{\mathbb R}^+$. Since $v\in\operatorname{L}^1_g([a,b))$,
	 there exists $\delta\in{\mathbb R}^+$ such that $\int_E v\operatorname{d}|\mu_g|<\varepsilon$ if $|\mu_g|(E)<\delta$. Then, for $\operatorname{var}_g[x,y]=|\mu_g|([x,y))<\delta$, we have that
\[
|h(y)-h(x)|= \left|\int_{[x,y)}v\operatorname{d}\mu_g\right| \leqslant \int_E v\operatorname{d}|\mu_g| <\varepsilon,
\]
so $h$ is $g$-continuous.

Now, take a family $\{(\alpha_\xi,\beta_\xi)\}_{\xi\in\Xi^+}$ of pairwise disjoint subintervals of $A_g^+$, a family $\{(\alpha_\xi,\beta_\xi)\}_{\xi\in\Xi^-}$ of pairwise disjoint subintervals of $A_g^-$, and sets $E^+\subset D_g^+$, $E^-\subset D_g^-$ such that
\[
\sum_{\xi\in\Xi^+}[g(\beta_\xi)-g(\alpha_\xi^+)] + \sum_{\xi\in\Xi^-}[g(\alpha_\xi^+)-g(\beta_\xi)]
+ \sum_{t\in E^+}[g(t^+)-g(t)] + \sum_{t\in E^-}[g(t)-g(t^+)] <\delta.
\]
Let us denote $B^+=\bigcup_{\xi\in\Xi^+}(\alpha_\xi,\beta_\xi)$ and $B^-=\bigcup_{\xi\in\Xi^-}(\alpha_\xi,\beta_\xi)$. { Observe that the sets $B^+$, $B^-$, $E^+$ and $E^-$ are pair-wise disjoint, so
\[|\mu_g|\left( B^+\cup B^- \cup E^+ \cup E^-\right)  = |\mu_g|(B^+) + |\mu_g|(B^-) + |\mu_g|(E^+) + |\mu_g|(E^-).\]
Furthermore, since $g$ is nondecreasing on the connected components of $B^+$, we have that $|\mu_g|(B^+)=\mu_g^+(B^+)$. By similar observations we get that
\[|\mu_g|\left( B^+\cup B^- \cup E^+ \cup E^-\right)  = \mu_g^+(B^+) + \mu_g^-(B^-) + \mu_g^+(E^+) + \mu_g^-(E^-).\]
Because  $g$ is left continuous, the measure in each of the intervals $(\alpha_\xi,\beta_\xi)$ is given by $|g(\beta_\xi)-g (\alpha_\xi^+)|$, so we have that}
\begin{align*}
&|\mu_g|\left( B^+\cup B^- \cup E^+ \cup E^-\right)  = \mu_g^+(B^+) + \mu_g^-(B^-) + \mu_g^+(E^+) + \mu_g^-(E^-)\\
 = & \sum_{\xi\in\Xi^+}[g(\beta_\xi)-g (\alpha_\xi^+)] -\sum_{\xi\in\Xi^-}[g(\beta_\xi)-g (\alpha_\xi^+)]
  + \sum_{t\in E^+}[g(t^+)-g (t)] - \sum_{t\in E^-}[g(t^+)-g (t)]
 <\delta.
 \end{align*}
 Thus,
 \[
 \sum_{\xi\in\Xi^+\cup \Xi^-}|h(\beta_\xi)-h(\alpha_\xi^+)| + \sum_{t\in E^+\cup E^-}|h(t^+)-h(t)|\\
 = \int_{B^+\cup B^-}v\operatorname{d}|\mu_g| + \int_{E^+\cup E^-}v\operatorname{d}|\mu_g| < \varepsilon.
\]
 Hence, $h\in\operatorname{AC}_g(I,{\mathbb R})$.
 \end{proof}

\begin{thm}\label{thmder}
Let $v\in\operatorname{L}^1_{g}([a,b))$ and $h(x):=\int_{[a,x)}v\operatorname{d}\mu_g$ \ for $x\in I$. Then, $h'_g=v$ \ $|\mu_g|$-a.e. on~$I$.
\end{thm}

\begin{proof}
 For the proof, we use~\cite[Theorem 2.4]{PoRo}, a theorem with the same statement but for a nondecreasing $g$. By definition, $A_g^+ = \bigcup_{\lambda\in\Lambda^+}I_{\lambda}$ and $g|_{I_{\lambda}}$ is nondecreasing for every $\lambda \in \Lambda^+$. {If we consider $h_{\lambda}(x)=\int_{[a_\lambda,x)}v\operatorname{d}\mu_{\widetilde g}$, where $\widetilde g(s)=g(s)$ for $s \in (a_\lambda,b]$ and $\widetilde g(a_\lambda):=g(a_\lambda^+)$, we have}
 \[
 h_{\lambda}(x)=\int_{{[}a_\lambda,x)}v\operatorname{d}\mu_{{\widetilde g}} =\int_{{[}a_\lambda,x)}v\operatorname{d}|\mu_{{\widetilde g}}|.
 \]
 It follows from~\cite[Theorem 2.4]{PoRo} that
 {$h_g'=(h_\lambda)'_g=(h_\lambda)'_{\widetilde g}=v$ holds} $|\mu_g|$-a.e. on $I_\lambda$.
 Since $\Lambda^+$ is countable, we deduce that
 $$
 h_g'=v \quad |\mu_g|\text{-a.e. on $A_g^+$}.
 $$

 Analogously, $A_g^- = \bigcup_{\lambda\in\Lambda^-}I_{\lambda}$ and $-g|_{I_{\lambda}}$ is nondecreasing for every $\lambda \in \Lambda^-$. So, for $x \in I_\lambda$,
 \[
 h_\lambda(x) = \int_{{[}a_\lambda,x)}v\operatorname{d}(-\mu_{\widetilde g}^-) = -\int_{{[}a_\lambda,x)}v\operatorname{d}\mu_{-{\widetilde g}},
 \]
 and
 $$
 h_g'=v \quad |\mu_g|\text{-a.e. on $A_g^-$.}
 $$

Let $t\in D_g$. Then
 \[
 h'_g(x)= \frac{h(x^+)-h(x)}{g(x^+)-g(x)} = \frac{\int_{[x,x^+)}v\operatorname{d}\mu_g}{g(x^+)-g(x)} = \frac{v(x)(g(x^+)-g(x))}{g(x^+)-g(x)} = v(x).
 \]
Hence, $h_g'=v$ in $D_g$.

Finally, $h_g' = v$ $|\mu_g|$-a.e. on $I$ since $I\backslash (A_g^+\cup A_g^- \cup D_g) = F_g\backslash D_g$ and $|\mu_g|(F_g\backslash D_g) = 0$.
\end{proof}

 Now, we show that $g$-absolutely continuous maps have bounded variation.

\begin{thm}\label{acbv}
If $h\in\operatorname{AC}_g(I,{\mathbb R})$, then $h\in\operatorname{BV}_-(I,{\mathbb R})$.
\end{thm}

\begin{proof}
By definition, $A_g^+ = \bigcup_{\lambda\in\Lambda^+}I_{\lambda}$, $A_g^- = \bigcup_{\lambda\in\Lambda^-}I_{\lambda}$ with $I_\lambda = (a_\lambda,b_\lambda)$. {Observe that, in general, for $\alpha\le \beta\le\gamma$, $\operatorname{var}_h[\alpha,\gamma]=\operatorname{var}_h[\alpha,\beta]+\operatorname{var}_h[\beta,\gamma]$. Hence, in order to bound the variation on $I$, we realize that it has to be less or equal than the sum of the variation in the intervals $(a_\lambda,b_\lambda]$ for $\lambda\in\Lambda^+\cup \Lambda^-$ plus the jumps at the points of $D_g$ plus an extra term where we make up for the fact that we have not considered the closed sets $[a_\lambda,b_\lambda]$, but omitted the point $a_\lambda$. Taking this into account, we have that}
\begin{equation}\label{eq:h-VB-0}
\begin{aligned}
\operatorname{var}_hI &\leqslant \sum_{\lambda\in\Lambda^+\cup \Lambda^-}\operatorname{var}_h (a_\lambda,b_\lambda] +
\sum_{\lambda\in\Lambda^+\cup \Lambda^- }|h(a_\lambda^+)-h(a_\lambda)| + \sum_{t\in D_g}|h(t^+)-h(t)|\\
&\leqslant \sum_{\lambda\in\Lambda^+\cup \Lambda^-}\operatorname{var}_h (a_\lambda,b_\lambda] +\sum_{t\in D_h}|h(t^+)-h(t)|.
\end{aligned}
\end{equation}
For $\varepsilon=1$, let $\delta\in{\mathbb R}^+$ be given in Definition~\ref{ac2}.

First, let us prove that $\sum_{t\in D_h}|h(t^+)-h(t)|<\infty$. By Proposition~\ref{lc}, $D_h\subset D_g$. Since $D_g$ is at most countable, without loss of generality, one can write $D_g = \{t_n : n \in {\mathbb N}\}$. Since $g \in \operatorname{CV}_-(I,{\mathbb R})$, we can fix $N_0 \in {\mathbb N}$ such that $\sum_{n=N_0+1}^\infty |g(t_n^+)-g(t_n)| < \delta$. Therefore, $\sum_{n=N_0+1}^\infty |h(t_n^+)-h(t_n)| < 1$. So,
\begin{equation}\label{eq:h-VB-1}
\sum_{t\in D_h}|h(t^+)-h(t)| < 1 + \sum_{n=1}^{N_0} |h(t_n^+)-h(t_n)| < \infty.
\end{equation}

Now, let us prove that
\[
\sum_{\lambda\in\Lambda^+\cup \Lambda^-}\operatorname{var}_h (a_\lambda,b_\lambda] < \infty.
\]
The set $\Lambda^+\cup \Lambda^-$ is at most countable and, without lost of generality, it can be written $\{\lambda_n : n \in {\mathbb N}\}$. Since $g$ has bounded variation, there exists $N\in {\mathbb N}$ such that
\[
\sum_{n=N+1}^\infty |g(b_{\lambda_n}) - g(a_{\lambda_n})| < \delta.
\]
For every $n > N$ and every $\{x_0^n,\dots,x_{m_n}^n\} \in \mathcal P(a_{\lambda_n},b_{\lambda_n})$, using the fact that $g$ is monotonic on $(a_{\lambda_n},b_{\lambda_n}]$, we get
\[
\sum_{n=N+1}^\infty\sum_{i=1}^{m_n}|g(x_i^n)-g(x_{i-1}^n)| = \sum_{n=N+1}^\infty |g(x_{m_n}^n)-g(x_0^n)| \leqslant \sum_{n=N+1}^\infty |g(b_{\lambda_n})-g(a_{\lambda_n})| < \delta.
\]
Thus,
\[
\sum_{n=N+1}^\infty\sum_{i=1}^{m_n}|h(x_i^n)-h(x_{i-1}^n)| < 1,
\]
and hence
\begin{equation}\label{eq:h-VB-2}
\sum_{n=N+1}^\infty\operatorname{var}_h (a_{\lambda_n},b_{\lambda_n}] \leqslant 1.
\end{equation}

Therefore, for $n = 1,\dots, N$, since $g$ is monotonic on $(a_{\lambda_n},b_{\lambda_n}]$, arguing as in~\cite[Proposition~5.3]{PoRo}, we deduce that
\begin{equation}\label{eq:h-VB-3}
\operatorname{var}_h (a_{\lambda_n},b_{\lambda_n}] < \infty.
\end{equation}

Combining~\eqref{eq:h-VB-0}--\eqref{eq:h-VB-3}, we conclude that $\operatorname{var}_hI <\infty$.
Finally, since $h$ is $g$-continuous, it is left continuous. Therefore, $h \in \operatorname{BV}_-(I,{\mathbb R})$.
\end{proof}

\begin{rem}
For $h\in\operatorname{AC}_g(I,{\mathbb R})$, since $h \in \operatorname{BV}_-(I,{\mathbb R})$ by the previous theorem, we can define the measures $\mu_h$ and $|\mu_h|$ associated to $h$ as we did in the preliminaries.
\end{rem}

\begin{lem}\label{lemac2}
 Let $h :I \to {\mathbb R}$ be a $g$-absolutely continuous function and $\mu_h$ its associated measure. Then, $\mu_h$ is $\mu_g^+$-absolutely continuous on $A_g^+$, (that is, $\mu_g^+(E)=0$ implies that $\mu_h(E)=0$ for any $\mu_g^+$-measurable set $E\subset A_g^+$) and $\mu_h$ is $\mu_g^-$-absolutely continuous on $A_g^-$.
 \end{lem}

 \begin{proof}
 We first show that, for any $\varepsilon\in{\mathbb R}^+$, there exists $\delta\in{\mathbb R}^+$ such that, for any open set $V \subset A_g^+$ satisfying $|\mu_g^+(V)|<\delta$, we have that $|\mu_h(V)|\leqslant \varepsilon$.
Fix an open set $V \subset A_g^+$ such that $|\mu_g^+(V)|<\delta$. The set $V$ can be written as a countable pairwise disjoint union of the form
$$
V=\bigcup_{\zeta\in Z^+}(r_\zeta,s_\zeta), \quad \text{where $(r_\zeta,s_\zeta)\subset I_\lambda$ for some $\lambda\in\Lambda^+$}.
 $$
Take $\widetilde r_\zeta\in(r_\zeta,s_\zeta)$ for every $\zeta\in Z^+$. Hence,
 \[
 \sum_{\zeta\in Z^+}[g(s_\zeta)-g(\widetilde r_\zeta)] = \mu_g\left( \bigcup_{\zeta\in Z^+}[\widetilde r_\zeta,s_\zeta)\right)  \leqslant \mu_g^+(V) <\delta,
 \]
 and, arguing as in the proof of the previous theorem, we get
 \[
 \sum_{\zeta\in Z^+}\operatorname{var}_h([\widetilde r_\zeta,s_\zeta)) \leqslant \varepsilon.
 \]
 Letting $\widetilde r_\zeta \to r_\zeta$, we obtain
 \[
 |\mu_h|(V) = \sum_{\zeta\in Z^+}\operatorname{var}_h((r_\zeta,s_\zeta)) \leqslant \varepsilon.
 \]

 Now, since $\mu_g^+$ and $|\mu_h|$ are outer regular, for any $\mu_g^+$-measurable set $E\subset A_g^+$, there exist open sets $V_n\subset A_g^+$ such that $E\subset V_n$ for every $n\in{\mathbb N}$ and
 \[
 \lim_{n\to\infty}\mu_g^+(V_n)=\mu_g^+(E) \quad\text{and}\quad \lim_{n\to\infty}|\mu_h|(V_n)=|\mu_h|(E).
 \]
 Now, for $E \subset A_g^+$ such that $\mu_g^+(E)=0$, one has $\mu_g^+(V_n) \to 0$ and, by the first part of the proof, $\lim_{n\to\infty}|\mu_h|(V_n)=0$. Therefore, $|\mu_h|(E)=0$ and hence, $\mu_h(E)=0$.

 Similarly, it can be shown that $\mu_h$ is $\mu_g^-$-absolutely continuous on $A_g^-$.
 \end{proof}

 We are ready to establish a Fundamental Theorem of Calculus for $g$-absolutely continuous maps.

\begin{thm}[Fundamental Theorem of Calculus for the Lebesgue-Stieltjes integral]
\label{ftc}
Let $g\in\operatorname{CV}_-(I,{\mathbb R})$ and $h: I\to{\mathbb R}$. The following statements are equivalent:
 \begin{enumerate}
 \item[(1)] $h\in\operatorname{AC}_g(I,{\mathbb R})$;
 \item[(2)] $h$ satisfies the following:
 \begin{enumerate}
 \item[(a)] there exists $h'_g$ $|\mu_g|$-a.e. in $I$,
 \item[(b)] $h'_g\in\operatorname{L}^1_g(I,{\mathbb R})$,
 \item[(c)] for each $t\in I$, we have that
 \[
 h(t)=h(a)+\int_{[a,t)}h'_g\operatorname{d}\mu_g.
 \]
 \end{enumerate}
 \end{enumerate}
\end{thm}

\begin{proof} $(1)\Rightarrow (2)$: We know from Lemma~\ref{lemac2} that $\mu_h$ is $\mu_g^+$-absolutely continuous on $A_g^+$. Hence, by the Radon-Nikodym Theorem~\cite[Theorem 5.3.2]{Ben}, there exists a unique function $v^+\in\operatorname{L}^1_{\mu_g^+}(A_g^+)$ such that
$$
\mu_h(E) = \int_E v^+\operatorname{d}\mu_g^+ \quad\text{for any $\mu_g^+$-measurable set $E\subset A_g^+$.}
$$
Similarly, there exists a unique function $v^-\in\operatorname{L}^1_{\mu_g^-}(A_g^-)$ such that
$$
\mu_h(E) = \int_E v^-\operatorname{d}\mu_g^- \quad\text{for any $\mu_g^-$-measurable set $E\subset A_g^-$.}
$$

Let us define $v : I \to {\mathbb R}$ by
$$
v(t) = \begin{dcases}
v^+(t) &\text{if $t \in A_g^+$,}\\
-v^-(t) &\text{if $t \in A_g^-$,}\\
\frac{h(t^+)-h(t)}{g(t^+)-g(t)} &\text{if $t \in D_g$,}\\
0 &\text{if $t \in F_g\backslash D_g$.}
\end{dcases}
$$
It follows from Theorem~\ref{acbv} that $v \in L^1_g(I,{\mathbb R})$ since
\begin{align*}
\int_I |v(t)|\operatorname{d}|\mu_g| &= \int_{A_g^+} |v^+(t)|\operatorname{d}\mu_g^+ + \int_{A_g^-} |v^-(t)|\operatorname{d}\mu_g^- + \int_{D_g} |v(t)|\operatorname{d}|\mu_g| + \int_{F_g\backslash D_g} |v(t)|\operatorname{d}|\mu_g| \\
&= \int_{A_g^+} |v^+(t)|\operatorname{d}\mu_g^+ + \int_{A_g^-} |v^-(t)|\operatorname{d}\mu_g^- + \sum_{t \in D_g}|h(t^+)-h(t)|\\
&< \infty.
\end{align*}
Moreover, for every $E \subset I$ $\mu_g$-measurable,
\begin{align*}
\mu_h(E) &= \mu_h(E\cap A_g^+) + \mu_h(E\cap A_g^-) + \mu_h(E \cap D_g) + \mu_h(E\cap(F_g\backslash D_g))\\
&=\int_{E\cap A_g^+} v^+(t)\operatorname{d}\mu_g^+ - \int_{E\cap A_g^-} -v^-(t)\operatorname{d}\mu_g^- + \sum_{t \in E\cap D_g} (h(t^+)-h(t))\\
&= \int_{E\cap A_g^+} v(t)\operatorname{d}\mu_g + \int_{E\cap A_g^-} v(t)\operatorname{d}\mu_g + \int_{E\cap D_g}v(t)\operatorname{d}\mu_g\\
&= \int_E v(t)\operatorname{d}\mu_g.
\end{align*}
In particular, for $E=[a,x)$, with $x\in[a,b]$ we have that
 \[
 h(x)-h(a)=\mu_h([a,x))=\int_{[a,x)}v\operatorname{d}\mu_g.
 \]
Finally, by Theorem~\ref{thmder}, $h'_g=v$ $|\mu_g|$-a.e. on $I$.

 $(2)\Rightarrow(1)$: That is a direct consequence of Theorem~\ref{thmint}.
 \end{proof}

 Here is another property of $g$-absolutely continuous functions.

 \begin{lem}\label{lemlip}
Let $h\in\operatorname{AC}_g(I,{\mathbb R})$ and $\varphi:h(I)\to{\mathbb R}$ be a Lipschitz function. Then $\varphi\circ h\in\operatorname{AC}_g(I,{\mathbb R})$.
\end{lem}

\begin{proof} Since $\varphi$ is Lipschitz, there exists $L\in{\mathbb R}^+$ such that $|\varphi(y)-\varphi(x)|\leqslant L|y-x|$ for every $x,y\in h(I)$. Fix $\varepsilon\in{\mathbb R}^+$ and consider a family $\{(\alpha_\xi,\beta_\xi)\}_{\xi\in\Xi^+}$ of pairwise disjoint subintervals of $A_g^+$, a family $\{(\alpha_\xi,\beta_\xi)\}_{\xi\in\Xi^-}$ of pairwise disjoint subintervals of $A_g^-$, and some sets $E^+\subset D_g^+$, $E^-\subset D_g^-$. Since $h\in\operatorname{AC}_g(I,{\mathbb R})$, $h$ is $g$-continuous and there exists $\delta\in{\mathbb R}^+$ such that if
\[
\sum_{\xi\in\Xi^+}[g(\beta_\xi)-g(\alpha_\xi^+)] + \sum_{\xi\in\Xi^-}[g(\alpha_\xi^+)-g(\beta_\xi)]
+ \sum_{t\in E^+}[g(t^+)-g(t)] + \sum_{t\in E^-}[g(t)-g(t^+)] <\delta,
\]
 then
 \[
 \sum_{\xi\in\Xi^+ \cup \Xi^-}|h(\beta_\xi)-h(\alpha_\xi^+)| + \sum_{t\in E^+\cup E^-}|h(t^+)-h(t)| <\frac{\varepsilon}{L}.
 \]
Thus,
\begin{align*}
& \sum_{\xi\in\Xi^+\cup \Xi^-}|\varphi\circ h(\beta_\xi)-\varphi\circ h(\alpha_\xi^+)| + \sum_{t\in E^+\cup E^-}|\varphi\circ h(t^+)-\varphi\circ h(t)| \\
\leqslant & L\left(\sum_{\xi\in\Xi^+\cup \Xi^-}|h(\beta_\xi)-h(\alpha_\xi^+)| + \sum_{t\in E^+\cup E^-}|h(t^+)-h(t)|\right) <\varepsilon.
\end{align*}
\end{proof}

We give an example of a $g$-absolutely continuous map with $g$ non monotonic which can be obtained from a solution of a system of differential equations with nondecreasing derivators.

\begin{exa} With this example, we show that a $g$-differential equation with a non monotonic $g$ can be obtained from systems of differential equations with nondecreasing derivators.

Let $p, q:I\to{\mathbb R}$ be left continuous and nondecreasing. Assume $p$ is $q$-absolutely continuous. Consider a system of differential equations of the form
 \begin{equation}\label{eqexa1}
 \begin{aligned}
 x'_{p}(t)= & f_1(t,x(t),y(t)),\\
 y'_{q}(t)= & f_2(t,x(t),y(t)),
 \end{aligned}
 \end{equation}
 for $t\in I$ and suppose it has a solution $(u,v)$. By assumption,
 $u$ is a solution of the differential equation:
 \begin{equation}\label{eqexa2}
 x'_{p}(t)= f_1(t,x(t),v(t)),
 \end{equation}
 and $u$ is $p$-absolutely continuous. That implies that $D_u\subset D_{p}\subset D_{q}$. Also, since $v$ is $q$-absolutely continuous, $D_v\subset D_{q}$. Let us assume that $D_v= D_{p}= D_{q}$, $v\in\operatorname{CV}_-(I,{\mathbb R})$, and $v$ and $v'_{q}$ are nonzero $|\mu_v|$-a.e. Let us point out that in general $v$ is not monotone.
 Now, we can consider the following $v$-differential equation:
 \begin{equation}\label{eqexa3}
 x'_v(t)=\frac{f_1(t,x(t),v(t))}{f_2(t,x(t),v(t))}{p}_{q}'(t).
 \end{equation}
 We have that, for $|\mu_v|$-a.e. $t\in I\backslash D_{v}$,
 \begin{align*}
 u_v'(t)= & \lim_{s\to t}\frac{u(s)-u(t)}{v(s)-v(t)}=
 \lim_{s\to t}\frac{u(s)-u(t)}{p(s)-p(t)} \frac{q(s)-q(t)}{v(s)-v(t)} \frac{p(s)-p(t)}{q(s)-q(t)}=
 \frac{u'_{p}(t)}{v'_{q}(t)}{p}_{q}'(t) \\
 = & \frac{f_1(t,u(t),v(t))}{f_2(t,u(t),v(t))}{p}_{q}'(t),
 \end{align*}
 and, for $t\in D_{v}$,
 \begin{align*}
 u_v'(t)= & \frac{u(t^+)-u(t)}{v(t^+)-v(t)}=
 \frac{u(t^+)-u(t)}{p(t^+)-p(t)} \frac{q(t^+)-q(t)}{v(t)-v(t^+)} \frac{p(t^+)-p(t)}{q(t^+)-q(t)}= \frac{u'_{p}(t)}{v'_{q}(t)}{p}_{q}'(t) \\
 = & \frac{f_1(t,u(t),v(t))}{f_2(t,u(t),v(t))}{p}_{q}'(t).
 \end{align*}
 That is, $u$ is a solution of~\eqref{eqexa3} as well.
 \end{exa}

\section{The exponential map}

{The computation of the exponential map for Stieltjes differential equations with positive derivator is a feat that was achieved in \cite{FP2016}. Here we generalize the construction of the exponential map for the case of a derivator $g$ which is allowed to change sign.}

Let $g \in\operatorname{CV}_-(I,{\mathbb R})$. In this section, we study the initial value problem for the linear $g$-differential equation:
\begin{equation}\label{eq:pb-lin}
\begin{aligned}
x_g'(t) &= c(t)x(t) \quad |\mu_g|\text{-a.e. $t \in I$,}\\
x(a) &= 1.
\end{aligned}
\end{equation}
To this aim, we introduce an exponential map of $c$ associated to $g$.

\begin{dfn}\label{dfnex} Let $c\in\operatorname{L}^1_g([a,b))$ be such that
\begin{equation}\label{ec1}
1+c(t)[g(t^+)-g(t)]>0,\ t\in[a,b)\cap D_g\ \text{ and }
\sum_{t\in[a,b)\cap D_g}\big|\ln(1+c(t)[g(t^+)-g(t)])\big|<\infty.
\end{equation}
 We define the {\em $g$-exponential of $c$ on $I$}, $e_{g,c,a}:I\to{\mathbb R}^+$, by
 \[
 e_{g,c,a}(t):= e^{\int_{[a,t)}\widetilde c(s)\operatorname{d}\mu_g},
 \]
 where
 $$
 \widetilde c(t):=\begin{dcases}
 c(t) &\text{if $t\in I\backslash D_g$,}\\
 \frac{\ln\left( 1+c(t)(g(t^+)-g(t))\right) }{g(t^+)-g(t)} &\text{if $t\in I\cap D_g$.}
 \end{dcases}
 $$
\end{dfn}

The proof of the following lemma is essentially the same as in~\cite[Lemma 6.2]{FP2016}, changing $g(s^+)-g(s)$ by $|g(s^+)-g(s)|$ since $g$ is non monotonic now.

\begin{lem}\label{lemct} The map $\widetilde c$ in Definition~\ref{dfnex} is such that $\widetilde c\in\operatorname{L}^1_g([a,b))$. So $e_{g,c,a}$ is well defined.
 \end{lem}

\begin{proof}
\begin{align*}\int_{[a,b)}|\widetilde c(s)|\operatorname{d} |\mu_g|= & \int_{[a,b)\backslash D_g}|\widetilde c(s)|\operatorname{d} \mu_g+\int_{[a,b)\cap D_g}|\widetilde c(s)|\operatorname{d} |\mu_g| \\= & \int_{[a,b)\backslash D_g}| c(s)|\operatorname{d} |\mu_g|+\sum_{s\in[a,b)\cap D_g}|\widetilde c(s)||g(s^+)-g(s)|\\ = & \int_{[a,b)\backslash D_g}| c(s)|\operatorname{d} |\mu_g|+\sum_{s\in[a,b)\cap D_g}\big|\ln\left( 1+c(s)[g(s^+)-g(s)]\right) \big|<\infty,
\end{align*}
because $c\in\operatorname{L}^1_g([a,b))$ and $\sum_{t\in[a,b)\cap D_g}\big|\ln(1+c(t)[g(t^+)-g(t)])\big|<\infty$.
\end{proof}

The $g$-exponential gives us a solution of~\eqref{eq:pb-lin}.

\begin{thm}\label{thmewd} Let $c\in\operatorname{L}^1_g([a,b))$. Then, $e_{g,c,a}$ solves the initial value problem~\eqref{eq:pb-lin}. Furthermore,
 \[
 e_{g,c,a}(t) = 1+\int_{[a,t)}c(s)e_{g,c,a}(s)\operatorname{d}\mu_g \quad\text{for every } t\in I.
 \]
\end{thm}

\begin{proof} By Lemma~\ref{lemct}, $\widetilde c\in\operatorname{L}^1_g([a,b))$. Therefore, by
Theorem~\ref{thmder}, for $h(t)=\int_{[a,t)}\widetilde c(s)\operatorname{d}\mu_g$, one has $h'_g=\widetilde c$ $|\mu_g|$-a.e. on~$I$. Hence, by Theorem~\ref{ftc}, $h\in\operatorname{AC}_g(I,{\mathbb R})$. Lemma~\ref{lemlip} implies that we have that $e_{g,c,a}\in\operatorname{AC}_g(I,{\mathbb R})$. Furthermore, by Corollary~\ref{ccr}, we can apply the chain rule and we get
$$
(e_{g,c,a})'_g(t)=e_{g,c,a}(t)h'_g(t)=e_{g,c,a}(t)\widetilde{c}(t) = e_{g,c,a}(t)c(t) \quad\text{$|\mu_g|$-a.e. on~$I\backslash D_g$.}
$$
Also, for $t \in D_g$, it is easy to verify that
\[
(e_{g,c,a})'_g(t) = \frac{e_{g,c,a}(t^+)-e_{g,c,a}(t)}{g(t^+)-g(t)} = e_{g,c,a}(t)c(t).
\]
 Now,
\[
e_{g,c,a}(a)=e^{\int_{[a,a)}\widetilde c(s)\operatorname{d}\mu_g}=1,
\]
and, since $e_{g,c,a}\in\operatorname{AC}_g(I,{\mathbb R})$, by Theorem~\ref{ftc},
\[
 e_{g,c,a}(t) = 1+\int_{[a,t)}c(s)e_{g,c,a}(s)\operatorname{d}\mu_g \quad\text{for every } t\in I.
 \]
\end{proof}

Now that we have covered the case~\eqref{ec1}, we can move to a more general one. In order to do that, consider $T_c^-:=\{t\in[a,b)\cap D_g\ :\ 1+c(t)[g(t^+)-g(t)]< 0\}$. Since $c\in\operatorname{L}^1_g([a,b))$, it can be shown that the set $T_c^-$ is finite (see~\cite[Lemma~6.4]{FP2016}), so assume $T_c^-=\{t_1,\dots,t_m\}$.

\begin{dfn}\label{dfnex2} Let $c\in\operatorname{L}^1_g([a,b))$ be such that
\begin{equation}\label{ec2}
1+c(t)[g(t^+)-g(t)]\ne 0,\ t\in[a,b)\cap D_g\ \text{ and }
\sum_{t\in[a,b)\cap D_g}\big|\ln|1+c(t)[g(t^+)-g(t)]|\big|<\infty.
\end{equation}
 We define the {\em $g$-exponential of $c$ on $I$}, $e_{g,c,a}:I\to{\mathbb R}^+$, as
 $$
 e_{g,c,a}(t):= \begin{dcases}e^{\int_{[a,t)}\widetilde c(s)\operatorname{d}\mu_g} &\text{if $t\in[a,t_1]$,} \\ (-1)^ke^{\int_{[a,t)}\widetilde c(s)\operatorname{d}\mu_g} &\text{if $t\in(t_k,t_{k+1}],\ k=1,\dots, m$,}
 \end{dcases}
 $$
 where
\begin{equation}\label{eqc}
 \widetilde c(t):=\begin{dcases}
 c(t) &\text{if $t\in I\backslash D_g$,}\\
 \frac{\ln|1+c(t)(g(t^+)-g(t))|}{g(t^+)-g(t)} &\text{if $t\in I\cap D_g$.}
 \end{dcases}
\end{equation}
\end{dfn}

\begin{rem} We use the same notation as in Definition~\ref{dfnex}, particularly $ e_{g,c,a}$ and $\widetilde c$, because the functions in Definitions~\ref{dfnex} and~\ref{dfnex2} coincide when~\ref{ec1} holds.
\end{rem}

Observe that Theorem~\ref{thmewd} also holds when condition~\eqref{ec2} holds as a direct consequence of Theorem~\ref{thmewd} holding for condition~\eqref{ec1} and the fact of $T_c^-$ being finite.

Finally, we can take this generalization one step further.

\begin{dfn}\label{dfnex3} Let $c\in\operatorname{L}^1_g([a,b))$ be such that
\begin{equation}\label{ec3}
\sum_{t\in[a,t_0)\cap D_g}\big|\ln|1+c(t)[g(t^+)-g(t)]|\big|<\infty,
\end{equation}
where $t_0=a$ if $T_c^0:=\{t\in[a,b)\cap D_g\ :\ 1 +c(t)[g(t^+)-g(t)]= 0\}=\emptyset$ and $t_0=\inf T_c^0$ otherwise.
 We define the {\em $g$-exponential of $c$ on $I$}, $e_{g,c,a}:I\to{\mathbb R}^+$, by
 \[
 e_{g,c,a}(t):= \begin{dcases}e^{\int_{[a,t)}\widetilde c(s)\operatorname{d}\mu_g}, & t\in[a,t_1]\cap[a,t_0], \\ (-1)^ke^{\int_{[a,t)}\widetilde c(s)\operatorname{d}\mu_g}, & t\in(t_k,t_{k+1}]\cap[a,t_0],\ k=1,\dots, m,\\ 0, & t\in(t_0,b], \end{dcases}
 \]
 where $\widetilde c$ is defined as in~\eqref{eqc}.
\end{dfn}

Observe that Theorem~\ref{thmewd} also holds when condition~\eqref{ec3} holds. The set $T_c^0$ is also finite, so we could write $t_0=\min T_c^0$.

\section{Existence results}

In this section, we present an existence result for the following system of $g$-differential equations:
\begin{equation}\label{gdifeq}
\begin{aligned}
x'_g(t) &=f(t,x(t)) \quad \mu_g\text{-a.e.}\ t\in[t_0,t_0+\tau],\\
 x(0) &=x_0,
 \end{aligned}
\end{equation}
where $I = [t_0,t_0+T]$ and $g=(g_1,\dots,g_n) :I\to{\mathbb R}^n$ is such that $g_j\in CV_-(I,{\mathbb R})$ for every $j=1,\dots,n$.

Similarly to the definition given in~\cite{FP2016}, a notion of $g$-Carath\'eodory function can be introduced in our context.

\begin{dfn}
Let $X$ be a nonempty subset of ${\mathbb R}^n$. We say that $f = (f_1,\dots,f_m):I\times X\to{\mathbb R}^m$ is \hbox{\emph{$g$-Carath\'eodory}} if it satisfies the following conditions for every $i=1,\dots,m$:
\begin{enumerate}
\item[(i)] for every $x\in X$, $f_i(\cdot,x)$ is $\mu_{g_i}$-measurable;
\item[(ii)] for $\mu_{g_i}$-a.e. $t\in I$, $f_i(t,\cdot)$ is continuous on $X$;
\item[(iii)] for every $r>0$, there exists $h_{i,r}\in\operatorname{L}^1_{g_i}([t_0,t_0+T),[0,\infty))$ such that
		$$
|f_i(t,x)|\leqslant h_{i,r}(t)\quad \text{for $\mu_{g_i}$-a.e. $t\in[t_0,t_0+T)$, for all $x\in X$ with $\|x\|\leqslant r$.}
$$
\end{enumerate}
\end{dfn}

We define $BC_g(I,{\mathbb R}^n):=\prod_{k=1}^nBC_{g_k}(I,{\mathbb R})$ with the supremum norm (let us recall that $BC_{g_k}(I,{\mathbb R})=BC((I,\tau_{g_k}),{\mathbb R})$, where $\tau_{g_k}$ is the topology generated by $g_k$). Arguing as in~\cite[Lemma~7.2]{FP2016}, it can be shown that the composition $f(\cdot,x(\cdot))$ is in $\operatorname{L}^1_g([t_0,t_0+T),{\mathbb R}^m) := \prod_{i=k}^m\operatorname{L}^1_{g_k}([t_0,t_0+T),{\mathbb R})$ for every $x \in BC_g(I,{\mathbb R}^n)$.

\begin{lem}\label{f-o-x} Let $X$ be a nonempty subset of ${\mathbb R}^n$ and $f:I\times X\to{\mathbb R}^m$ a $g$-Carath\'eodory function. Then, for every $x \in BC_g(I,{\mathbb R}^n)$, the map $f(\cdot,x(\cdot))$ is in $\operatorname{L}^1_g([t_0,t_0+T),{\mathbb R}^m)$.
\end{lem}

Now, we can establish the existence of $x=(x_1,\dots,x_n)$ a $g$-absolutely continuous solution of~\eqref{gdifeq}, i.e. $x_i$ is $g_i$-absolutely continuous for every $i=1,\dots,n$. The proof of the next theorem is analogous to that of~\cite[Theorem~4.5]{PoMa2}. We present it here in detail for sake of completeness.

\begin{thm}\label{exres2}
 Let $r\in{\mathbb R}^+$ and $f:[t_0,t_0+T]\times\overline{B(x_0,r)}\to{\mathbb R}^n$ a $g$-Carath\'eodory function. Then there exists $\tau\in(0,T]$ such that~\eqref{gdifeq}
 has a $g$-absolutely continuous solution defined on $[t_0,t_0+\tau]$.
\end{thm}

\begin{proof} Let $R:=r+\|x_0\|$. Since $f$ is $g$-Carath\'eodory, we have that, for every $j=1,\dots,n$, there exists a function $h_j\in\operatorname{L}^1_{g_j}([t_0,t_0+T),[0,\infty))$ such that
$|f_j(t,x)| \leqslant h_j(t)$, for $\mu_{g_j}$-a.e. $t\in[t_0,t_0+T)$ and all $x\in X$ with $\|x\|\leqslant r$. Fix $\tau\in(0,T]$ such that
\[
\max_{j=1,\dots,n}\int_{[t_0,t_0+\tau)} h_j\operatorname{d} |\mu_{g_j}|\leqslant r.
\]
Let
\[
X:=\{u\in BC_g([t_0,t_0+\tau),{\mathbb R}^n)\ : \|u-x_0\|_\infty\leqslant r\}.
\]
The set $X$ is a nonempty closed convex subset of $BC_g([t_0,t_0+\tau),{\mathbb R}^n)$.
For $u\in BC_g([t_0,t_0+\tau),{\mathbb R}^n)$ define $F(u) = (F_1(u),\dots,F_n(u))$ with
$$
F_j(u)=x_{j,0} +\int_{[t_0,t_0+t)}f_j(s,u(s))\operatorname{d} \mu_{g_j}, \quad \text{for $j=1,\dots,n$}.
$$
The previous lemma implies that $F$ is well defined. It is easy to verify that every fixed point of $F$ is a solution of~\eqref{gdifeq} and vice-versa.

Now, we prove certain properties of $F$ that will allow us to use Schauder's Fixed Point Theorem.

\emph{\textbullet\ $F$ maps $X$ to $X$.}

Let us show that $Fu\in BC_g([t_0,t_0+\tau),{\mathbb R}^n)$ for $u\in BC_g([t_0,t_0+\tau),{\mathbb R}^n)$. Let $u\in BC_g([t_0,t_0+\tau),{\mathbb R}^n)$. By~\ref{f-o-x}, we have that $f(\cdot,u(\cdot))\in\operatorname{L}^1_g([t_0,t_0+\tau),{\mathbb R}^n)$. Thus, by Theorem~\ref{thmint}, we have that $Fu\in BC_g([t_0,t_0+\tau),{\mathbb R}^n)$.

Let us show that $Fu\in X$ for $u\in X$. Let $u\in X$. We have that $\|u-x_0\|\leqslant r$. Then, for any $j=1,\dots,n$ and $t\in[t_0,t_0+\tau)$,
\[
|(Fu-x_0)_j(t)|\leqslant\int_{[t_0,t)}|f_j(s,u(s))| \operatorname{d} |\mu_{g_j}|(s)\leqslant\int_{[t_0,t_0+\tau)} h_j\operatorname{d} |\mu_{g_j}|\leqslant r,
\]
so $Fu\in X$.

\emph{\textbullet\ $F$ is continuous.}

Let $(u_n)_{n\in{\mathbb N}}\subset X$, $u\in X$ be such that $u_n\to u$. Then $u_n(s)\to u_n(s)$ for every $s\in[t_0,t_0+\tau)$, and, since $f_j$ is $g_j$-Carath\'eodory, $f_j(s,u_n(s))\to f_j(s,u(s))$ $\mu_{g_j}$-a.e. Now, for any $j=1,\dots,n$ and $t\in[t_0,t_0+\tau)$,
\begin{align*}
|F_ju_n(t)-F_ju(t)|\leqslant\int_{[t_0,t_0+\tau)}|f_j(s,u_n(s))-f_j(s,u(s))| \operatorname{d} |\mu_{g_j}|(s)\leqslant 2r.
\end{align*}
Thus, by the Lebesgue Dominated Convergence Theorem,
\begin{align*}
\lim_{n\to\infty}\|Fu_n-Fu\|_\infty\leqslant & \lim_{n\to\infty} \int_{[t_0,t_0+\tau)}|f_j(s,u_n(s))-f_j(s,u(s))| \operatorname{d} |\mu_{g_j}| \\ \leqslant & \int_{[t_0,t_0+\tau)}\lim_{n\to\infty}|f_j(s,u_n(s))-f_j(s,u(s))| \operatorname{d} |\mu_{g_j}|=0
\end{align*}
because $f_j(s,u_n(s))\to f_j(s,u(s))$ $\mu_{g_j}$-a.e.

\emph{\textbullet\ $F$ is compact.}

Since $X$ is bounded $F(X)\subset X$ is uniformly bounded. Furthermore, $F(X)$ is uniformly $g$-e\-qui\-con\-tin\-u\-ous. Indeed, let $\varepsilon\in{\mathbb R}^+$ and $j=1,\dots,n$. Since $h_j\in\operatorname{L}^1_{g_j}([t_0,t_0+\tau),{\mathbb R})$, we have that there exists $\delta\in{\mathbb R}^+$ such that if $|\mu_{g_j}|(E)<\delta$, then $\int_{E}f_j(s,u(s))\operatorname{d} |\mu_{g_j}|(s)<\varepsilon$. Hence, if $t,s\in[t_0,t_0+\tau)$ are such that $\operatorname{var}_{g_j}[t,s]=|\mu_{g_j}|([t,s])<\delta$, then
\[
|F_ju(t)-F_ju(s)|=\left|\int_{[s,t)}f_j(s,u(s))\operatorname{d} \mu_{g_j}(s)\right|\leqslant \int_{[s,t)}h_j\operatorname{d} |\mu_{g_j}|<\varepsilon.
\]
So, $F(X)$ is uniformly $g$-equicontinuous. It follows from Corollary~\ref{cor:precompact} that $F(X)$ is precompact and hence, $F$ is a compact operator.

Hence, we are in the hypotheses of Schauder's Fixed Point Theorem and we can ensure that $F$ has a fixed point.
\end{proof}

Other results similar to Theorem~\ref{exres2} can be proven in an standard way. The reader may refer, for instance, to Theorems~4.2,~4.3,~4.4 and~4.8 in~\cite{PoMa2}.

\section{An application to fluid stratification on buoyant miscible jets and plumes}

In~\cite{Camassa}, the authors modeled fluid stratification on buoyant miscible jets and plumes using the classical model of Morton et al.~\cite{Morton}, which is given by
\begin{equation}\label{me}\begin{aligned}
(b^2w)'= & 2\alpha bw,\\
(b^2w^2)'= & 2g\lambda^2b^2\theta,\\
(b^2w\theta)'= & \rho'\frac{b^2w}{\Lambda\rho_b},
\end{aligned}
\end{equation}
where the independent variable, $z$, is the height, $b(z)$ the jet width, $w(z)$ the vertical jet velocity, $\theta(z)=(\rho(z)-\rho_j(z))/\rho_b$ the density anomaly, $\rho(z)$ the ambient background density, $g$ the acceleration due to gravity, $\alpha$ the entrainment coefficient, $\lambda$ the mixing coefficient and $\Lambda=\lambda^2(1+\lambda^2)$. With the change of variables $q=b^2w$, $m=b^4w^4$ and $\beta=b^2w\theta$ (volume, momentum and buoyancy fluxes),~\eqref{me} becomes
\begin{equation}\begin{aligned}
q'= & 2\alpha m^{\frac{1}{4}},\\
m'= & 4g\lambda^2 q\beta,\\
\beta'= & \rho'\frac{q}{\Lambda\rho_b}.
\end{aligned}
\end{equation}

In~\cite{Camassa}, the authors consider the particular case where $\rho$ suffers a sudden change at $z=L$ (a jump) and thus $\rho'$ can be approximated by a Dirac delta function at $L$ (see~\cite[Eq. 2.35c]{Camassa}). In our approach, we will not consider distributions in order to obtain jumps. We will just consider the $\rho$ derivative of $\beta$ which will codify any possible behavior of $\rho$. Hence, this will lead us to consider the system of Steiltjes differential equations
\begin{equation}\label{emf}\begin{aligned}
q'= & f_1(q,m,\beta)=A m^{\frac{1}{4}},\\
m'= & f_2(q,m,\beta)=Bq\beta,\\
\beta'_\rho= & f_3(q,m,\beta)=C q,
\end{aligned}
\end{equation}
where $A=2\alpha$, $B=4g\lambda^2$ and $C=1/(\Lambda\rho_b)$. This is a system of differential equations with $g_j$-derivatives where $g_1$ and $g_2$ are the identity and $g_3=\rho\in CV_-({\mathbb R},{\mathbb R}^+)$. Now fix $x_0=(q_0,m_0,\beta_0)\in({\mathbb R}^+)^3$. It is clear that there exists $r\in{\mathbb R}^+$ such that $f_j:[0,T]\times\overline{B(x_0,r)}\to{\mathbb R}$ is $g_j$-Carath\'eodory
(just take $r<m_0$). Thus, by Theorem~\ref{exres2}, there exists $\tau\in(0,T]$ such a
that~\eqref{emf} has a $g$-absolutely continuous solution defined on $[t_0,t_0+\tau]$.

\section*{Declarations}
\subsection*{Availability of data and material}
Not applicable.

\subsection*{Competing interests}
The authors declare that they have no competing interests.
\subsection*{Funding}
	Marl\`ene Frigon was partially supported by NSERC Canada.

F. Adrián F. Tojo was partially supported by Ministerio de Econom\'ia y Competitividad, Spain, and FEDER, project MTM2013-43014-P, and by the Agencia Estatal de Investigaci\'on (AEI) of Spain under grant MTM2016-75140-P, co-financed by the European Community fund FEDER.

\subsection*{Authors' contributions}

All authors contributed equally to the different parts of the manuscript. All authors read and approved the final manuscript.

\subsection*{Acknowledgements}

F. Adrián F. Tojo would like to acknowledge his gratitude towards Prof. Marlène Frigon and the \emph{Département de mathématiques et de statistique} of the \textit{Université de Montréal} for their warm accueil during his stay at the aforementioned University, time during which this article was written.


\end{document}